\documentclass[10pt,twoside,reqno]{amsart}
\usepackage{amssymb}
\usepackage{multirow}
\calclayout

\numberwithin{equation}{section}
\makeatletter

\renewcommand{\@secnumfont}{\bfseries}

\renewcommand{\section}{\@startsection{section}{1}%
  {0mm}{.7\linespacing\@plus\linespacing}{.5\linespacing}
  {\normalfont\bfseries\centering}}

\newcommand{\bibsection}{\@startsection{section}{1}%
  {0mm}{.7\linespacing\@plus\linespacing}{.5\linespacing}
  {\normalfont\scshape\centering}}

\renewcommand{\@biblabel}[1]{#1.}

\begin{document}

\vspace{1.3cm}

\title {Degenerate Eulerian numbers and polynomials}

\author{Taekyun Kim}
\address{Department of Mathematics, Kwangwoon University, Seoul 139-701, Republic
	of Korea}
\email{tkkim@kw.ac.kr}

\author{Dae San Kim}
\address{Department of Mathematics, Sogang University, Seoul 121-742, Republic of Korea}
\email{dskim@sogang.ac.kr}

\subjclass[2010]{11B83; 11S80}
\keywords{Degenerate Eulerian numbers and polynomials}

\begin{abstract} 
In this paper, we study the degenerate Eulerian polynomials and numbers and give some new and interesting identities associated with several special numbers and polynomials.
\end{abstract}
\maketitle

\section{Introduction}
In combinatorics, the Eulerian number $\genfrac<>{0pt}{}{n}{m}$, is the number of permutations of the numbers 1 to $n$ in which exactly $m$ elements are greater than the previous element.

Indeed, the generating function of Eulerian numbers is given by
\begin{equation}\begin{split}\label{01}
\left( \sum_{k=0}^\infty (k+1)^n x^k \right) (1-x)^{n+1} = \sum_{m=1}^\infty \genfrac<>{0pt}{}{n}{m-1} x^m,\quad (\textnormal{see} \,\, [7,10]).
\end{split}\end{equation}
Thus, by \eqref{01}, we get
\begin{equation}\begin{split}\label{02}
\genfrac<>{0pt}{}{n}{m} = \sum_{l=0}^{m+1} {n+1 \choose l} (-1)^l (m+1-l)^n, \,\,(n \in \mathbb{N}, m \geq 0).
\end{split}\end{equation}
From \eqref{02}, we note that
\begin{equation*}\begin{split}
\genfrac<>{0pt}{}{n}{m} &= \sum_{k=0}^{m+1} {n+1 \choose k} (-1)^k (m+1-k)^n \\
&= \sum_{k=0}^{m+1} {n+1 \choose k} (-1)^k (m+1-k)^{n-1} (m+1-k)\\
&= (m+1)\sum_{k=0}^{m+1} {n+1 \choose k} (-1)^k (m+1-k)^{n-1}\\
&\quad- \sum_{k=1}^{m+1} {n+1 \choose k}k (-1)^k (m+1-k)^{n-1} \qquad\qquad\qquad\qquad\qquad\qquad\qquad\qquad
\end{split}\end{equation*}
\begin{equation}\begin{split}\label{03}
&= (m+1) \sum_{k=0}^{m+1} {n+1 \choose k} (-1)^k (m+1-k)^{n-1}\\
&\quad - (n+1) \sum_{k=1}^{m+1} {n \choose k-1} (-1)^k (m+1-k)^{n-1}\\
&= (m+1) \sum_{k=0}^{m+1} {n+1 \choose k}(-1)^k (m+1-k)^{n-1}+(n+1) \sum_{k=0}^m {n \choose k} (-1)^k (m-k)^{n-1}\\
&=(m+1) \sum_{k=0}^{m+1} \left( {n \choose k} + {n \choose k-1} \right)(-1)^k (m+1-k)^{n-1} + (n+1) \genfrac<>{0pt}{}{n-1}{m-1}\\
&= (m+1) \genfrac<>{0pt}{}{n-1}{m} \\
&\quad+ (m+1) \sum_{k=1}^{m+1}{n \choose k-1}(-1)^k (m+1-k)^{n-1} + (n+1) \genfrac<>{0pt}{}{n-1}{m-1}\\
&= (m+1) \genfrac<>{0pt}{}{n-1}{m} +(m+1) \sum_{k=0}^m {n \choose k} (-1)^{k-1} (m-k)^{n-1} + (n+1)\genfrac<>{0pt}{}{n-1}{m-1}\\
&= (m+1) \genfrac<>{0pt}{}{n-1}{m} - (m+1) \genfrac<>{0pt}{}{n-1}{m-1} + (n+1) \genfrac<>{0pt}{}{n-1}{m-1}\\
&= (n-m)\genfrac<>{0pt}{}{n-1}{m-1} + (m+1)\genfrac<>{0pt}{}{n-1}{m}.
\end{split}\end{equation}

By \eqref{03}, we obtain the recurrence relation for Eulerian numbers as follows:
\begin{equation}\begin{split}\label{04}
\genfrac<>{0pt}{}{n}{m} = (n-m)\genfrac<>{0pt}{}{n-1}{m-1} + (m+1) \genfrac<>{0pt}{}{n-1}{m},\quad (\textnormal{see} \,\, [3,7,10]).
\end{split}\end{equation}
As is well known, the Eulerian polynomials, $A_n(t)$, $(n \geq 0)$, are defined by the generating function 
\begin{equation}\begin{split}\label{05}
\frac{1-t}{e^{x(t-1)}-t} = e^{A(t)x} = \sum_{n=0}^\infty A_n(t) \frac{x^n}{n!},
\end{split}\end{equation}
with the usual convention about replacing $A^n(t)$ by $A_n(t)$. From \eqref{04}, we note that
\begin{equation}\begin{split}\label{06}
\big( A(t) + (t-1) \big)^n - tA_n(t) = (1-t) \delta_{0,n}, \,\,(n \geq 0),
\end{split}\end{equation}
where $\delta_{n,k}$ is the Kronecker's symbol (see [7]).
From \eqref{03}, \eqref{04}, \eqref{05} and \eqref{06}, we note that
\begin{equation}\begin{split}\label{07}
A_n(t) = \sum_{l=0}^n \genfrac<>{0pt}{}{n}{l} t^l, \,\,(n \geq 0),\quad (\textnormal{see} \,\, [ 3,7,10]).
\end{split}\end{equation}
The first few Eulerian polynomials are given by
\begin{equation}\begin{split}\label{08}
&1+t+t^2+t^3+\cdots =\frac{1}{1-t} = \frac{A_0(t)}{1-t},\\
&1+2t+3t^2+4t^3+\cdots = \frac{1}{(1-t)^2} = \frac{A_1(t)}{(1-t)^2},\\
&1+2^2t+3^2t^2+4^2t^3+\cdots = \frac{1+t}{(1-t)^3} = \frac{A_2(t)}{(1-t)^3},\\
&1+2^3t+3^3t^2+4^3t^3+\cdots = \frac{1+4t+t^2}{(1-t)^4} = \frac{A_3(t)}{(1-t)^4}.
\end{split}\end{equation}
The Worpitzky's identity expresses $x^n$ as the linear combination of Eulerian numbers with binomial coefficients as follows:
\begin{equation}\begin{split}\label{09}
x^n = \sum_{k=0}^{n-1} \genfrac<>{0pt}{}{n}{k} {x+k \choose n},\quad (\textnormal{see} \,\, [ 3,4,5,6,7,9,10]).
\end{split}\end{equation}
From \eqref{06}, we note that
\begin{equation}\begin{split}\label{10}
A_0(t) =1, A_n(t) = \frac{1}{t-1} \sum_{l=0}^{n-1} {n \choose l} A_l(t) (t-1)^{n-l}, \,\,(n \geq 1),
\end{split}\end{equation}
and
\begin{equation}\begin{split}\label{11}
\sum_{k=1}^m k^m t^k = \sum_{i=1}^n (-1)^{n+i} {n \choose i} \frac{t^{m+1}A_{n-i}(t)}{(t-1)^{n-i+1}}m^i + (-1)^n \frac{t(t^m-1)}{(t-1)^{n+1}}A_n(t),
\end{split}\end{equation}
where $m \geq 1$ and $n \geq 0$ (see [7,10]).

In [6], the degenerate ordered Bell polynomials are defined by the generating function 
\begin{equation}\begin{split}\label{12}
\frac{1}{2-(1+\lambda t)^{\frac{1}{\lambda }}}(1+\lambda t)^{\frac{x}{\lambda }}= \sum_{n=0}^\infty b_{n,\lambda }(x) \frac{t^n}{n!}.
\end{split}\end{equation}
When $x=0$, $b_{n,\lambda }=b_{n,\lambda }(0)$ are called the degenerate ordered Bell numbers.
It is well known that the Frobenius-Euler polynomials are given by the generating function 
\begin{equation}\begin{split}\label{13}
\frac{1-u}{e^t-u} e^{xt} = \sum_{n=0}^\infty H_n(x|u) \frac{t^n}{n!}, 
\end{split}\end{equation}
where $u \neq 1$. (see [8]). When $x=0$, $H_n(u) = H_n(0|u)$ are called the Frobenius-Euler numbers. Recently, several authors have studied some interesting extensions and modifications of Eulerian polynomials and numbers (see [1-12]).

In this paper, we study the degenerate Eulerian polynomials and numbers, which are due to Carlitz (see [1]), and give some new and interesting identities for these numbers and polynomials associated with several special numbers and polynomials.

\section{Degenerate Eulerian polynomials and numbers}

We recall that the Stirling numbers of the first kind and of the second kind are defined by the generating function as follows:
\begin{equation}\begin{split}\label{14}
\frac{1}{k!} \big( \log(1+t)\big)^k = \sum_{n=k}^\infty S_1(n,k) \frac{t^n}{n!},
\end{split}\end{equation}
and 
\begin{equation}\begin{split}\label{15}
\frac{1}{k!} \big(e^t-1 \big)^k = \sum_{n=k}^\infty S_2(n,k) \frac{t^n}{n!},\quad (\textnormal{see} \,\, [11]).
\end{split}\end{equation}
For $\lambda \in \mathbb{R}$, we consider the degenerate Eulerian polynomials given by the generating function 
\begin{equation}\begin{split}\label{16}
\frac{1-t}{(1+\lambda x)^{\frac{t-1}{\lambda }}-t} = \sum_{n=0}^\infty A_{n,\lambda }(t) \frac{x^n}{n!}.
\end{split}\end{equation}
Note that $\lim_{\lambda  \rightarrow 0} A_{n,\lambda }(t) = A_n(t)$, $(n \geq 0)$. From \eqref{16}, we have
\begin{equation}\begin{split}\label{17}
1-t &= \left( \sum_{n=0}^\infty A_{n,\lambda }(t) \frac{x^n}{n!} \right)\cdot \left( (1+\lambda x)^{\frac{t-1}{\lambda }}-t \right)\\
&= \left( \sum_{k=0}^\infty A_{k,\lambda }(t) \frac{x^k}{k!} \right) \left( \sum_{m=0}^\infty {\frac{t-1}{\lambda } \choose m} \lambda ^m x^m-t \right)\\
&= \left( \sum_{k=0}^\infty A_{k,\lambda }(t) \frac{x^k}{k!} \right) \left( \sum_{m=0}^\infty (t-1)_{m,\lambda } \frac{x^m}{m!} -t \right)\\
&= \sum_{n=0}^\infty \left( \sum_{k=0}^n {n \choose k} A_{k,\lambda }(t) (t-1)_{n-k,\lambda }- tA_{n,\lambda }(t) \right) \frac{x^n}{n!},
\end{split}\end{equation}
where $(x)_{n,\lambda } = x(x-\lambda )\cdots (x-(n-1)\lambda )$, $(n \geq 1)$, $(x)_{0,\lambda }=1$.

Comparing the coefficients on both sides of \eqref{17}, we get
\begin{equation}\begin{split}\label{18}
\sum_{k=0}^n {n \choose k} A_{k,\lambda }(t) (t-1)_{n-k,\lambda }- tA_{n,\lambda }(t) = (1-t)\delta_{0,n}.
\end{split}\end{equation}
Thus, from \eqref{18}, we have
\begin{equation}\begin{split}\label{19}
\sum_{k=0}^{n-1} {n \choose k} A_{k,\lambda }(t) (t-1)_{n-k,\lambda } =(t-1) A_{n,\lambda }(t), \,\,(n \geq 1), \,\, A_{0,\lambda }(t) =1.
\end{split}\end{equation}
For $n \geq 1$, we have
\begin{equation}\begin{split}\label{20}
A_{n,\lambda }(t) = \frac{1}{t-1} \sum_{k=0}^{n-1} {n \choose k} A_{k,\lambda }(t) (t-1)_{n-k,\lambda }.
\end{split}\end{equation}
From \eqref{16}, we note that
\begin{equation}\begin{split}\label{21}
\sum_{n=0}^\infty A_{n,\lambda }(t) \frac{x^n}{n!} &=
\frac{1-t}{(1+\lambda x)^{\frac{t-1}{\lambda }}-t} = \frac{1-t}{ e^{\frac{t-1}{\lambda }\log(1+\lambda x)}   -t}\\
&= \sum_{k=0}^\infty A_k(t) \frac{\lambda ^{-k}}{k!} \big( \log(1+\lambda x) \big)^k\\
&= \sum_{k=0}^\infty A_k(t) \lambda ^{-k} \sum_{n=k}^\infty S_1(n,k) \frac{\lambda ^n x^n}{n!}\\
&= \sum_{n=0}^\infty \left( \sum_{k=0}^n A_k(t) \lambda ^{n-k} S_1(n,k) \right) \frac{x^n}{n!}.
\end{split}\end{equation}
Thus, by comparing the coefficients on both sides of \eqref{21}, we get
\begin{equation}\begin{split}\label{22}
A_{n,\lambda }(t) =  \sum_{k=0}^n A_k(t) \lambda ^{n-k} S_1(n,k),\,\,(n \geq 0).
\end{split}\end{equation}
In view of \eqref{07}, we define the degenerate Eulerian polynomials by
\begin{equation}\begin{split}\label{23}
A_{n,\lambda }(t) = \sum_{l=0}^n \genfrac<>{0pt}{}{n}{l}_\lambda  t^l.
\end{split}\end{equation}
Thus, we easily get $\lim_{\lambda  \rightarrow 0} \genfrac<>{0pt}{}{n}{l}_\lambda  = \genfrac<>{0pt}{}{n}{l}$, $(n \geq 0)$.
From \eqref{07}, \eqref{22} and \eqref{23}, we have
\begin{equation}\begin{split}\label{24}
\sum_{l=0}^n \genfrac<>{0pt}{}{n}{l}_\lambda  t^l &= A_{n,\lambda }(t) = \sum_{k=0}^n A_k(t) \lambda ^{n-k} S_1(n,k)\\
&= \sum_{k=0}^n \sum_{l=0}^k \genfrac<>{0pt}{}{k}{l} t^l \lambda ^{n-k} S_1(n,k)\\
&= \sum_{l=0}^n \left( \sum_{k=l}^n \genfrac<>{0pt}{}{k}{l} \lambda ^{n-k} S_1(n,k) \right)t^l.
\end{split}\end{equation}
Comparing the coefficients on both sides of \eqref{24}, we obtain
\begin{equation}\begin{split}\label{25}
\genfrac<>{0pt}{}{n}{l}_\lambda  = \sum_{k=l}^n \genfrac<>{0pt}{}{k}{l} \lambda ^{n-k} S_1(n,k),\,\,(0 \leq l \leq n).
\end{split}\end{equation}
By \eqref{12} and \eqref{16}, we get
\begin{equation}\begin{split}\label{26}
\sum_{n=0}^\infty b_{n,\lambda } \frac{x^n}{n!} = \frac{1}{2-(1+\lambda x)^{\frac{1}{\lambda }}} = \sum_{n=0}^\infty A_{n,\lambda }(2) \frac{x^n}{n!}.
\end{split}\end{equation}
Thus, by \eqref{24}, we have
\begin{equation}\begin{split}\label{27}
b_{n,\lambda } = A_{n,\lambda }(2) &= \sum_{l=0}^n  \genfrac<>{0pt}{}{n}{l}_\lambda 2^l \\
&= \sum_{l=0}^n \sum_{k=l}^n  \genfrac<>{0pt}{}{k}{l} \lambda ^{n-k} S_1(n,k) 2^l,
\end{split}\end{equation}
where $0 \leq l \leq n$.
For $0 \leq l \leq n$, we have
\begin{equation}\begin{split}\label{28}
 \genfrac<>{0pt}{}{n}{l}_\lambda  = \sum_{k=l}^n \sum_{m=0}^{l+1} {k+1 \choose m} (-1)^m (l+1-m)^k \lambda ^{n-k} S_1(n,k).
\end{split}\end{equation}
Note that
\begin{equation*}\begin{split}
\lim_{\lambda  \rightarrow 0}  \genfrac<>{0pt}{}{n}{l}_\lambda &=
 \sum_{m=0}^{l+1} {n+1 \choose m} (-1)^m (l+1-m)^n\\
 &=  \genfrac<>{0pt}{}{n}{l},\,\,(0 \leq l \leq n).
\end{split}\end{equation*}
From \eqref{21}, we can derive the following equation:
\begin{equation}\begin{split}\label{29}
\sum_{n=0}^\infty A_{n,\lambda }(t) \frac{x^n}{n!}&=
\frac{1-t}{(1+\lambda x)^{\frac{t-1}{\lambda }}-t} = \frac{1-t}{ e^{\frac{t-1}{\lambda }\log(1+\lambda x)}   -t}\\
&= \sum_{n=0}^\infty H_n(t) \frac{1}{n!} \left(\frac{t-1}{\lambda }\right)^n \big( \log(1+\lambda x)\big)^n \\
&= \sum_{k=0}^\infty H_k(t) \left( \frac{t-1}{\lambda } \right)^k \sum_{n=k}^\infty S_1(n,k) \frac{\lambda ^nx^n}{n!}\\
&= \sum_{n=0}^\infty \left( \sum_{k=0}^n (t-1)^k H_k(t) \lambda ^{n-k} S_1(n,k) \right) \frac{x^n}{n!},
\end{split}\end{equation}
where $H_n(t)$ is the Frobenius-Euler numbers. By \eqref{29}, we get
\begin{equation}\begin{split}\label{30}
A_{n,\lambda }(t) = \sum_{k=0}^n \lambda ^{n-k} S_1(n,k) H_k(t)(t-1)^k ,\,\,(n \geq 0).
\end{split}\end{equation}
Let us take $t=2$. Then we have
\begin{equation}\begin{split}\label{31}
b_{n,\lambda } = \sum_{k=0}^n \lambda ^{n-k} S_1(n,k) H_k(2),\,\,(n \geq 0).
\end{split}\end{equation}

\section{Further remark}

Let $p$ be an odd prime number. Throughout this section, $\mathbb{Z}_p$, $\mathbb{Q}_p$ and $\mathbb{C}_p$ will denote the ring of $p$-adic integers, the field of $p$-adic rational numbers and the completion of the algebraic closure of $\mathbb{Q}_p$, respectively. The $p$-adic norm is normalized so that $|p|_p = \frac{1}{p}$. Let $q$ be an indeterminate in $\mathbb{C}_p$ such that
$|1-q|_p < p^{-\frac{1}{p-1}}$. As notations, the $q$-numbers are defined by
\begin{equation*}\begin{split}
[x]_q = \frac{1-q^x}{1-q},\quad \text{and}\quad [x]_{-q} = \frac{1-(-q)^x}{1+q}.
\end{split}\end{equation*}
Let $f$ be a continuous function on $\mathbb{Z}_p$. Then the fermionic $p$-adic $q$-integral on $\mathbb{Z}_p$ is defined as
\begin{equation}\begin{split}\label{32}
I_{-q}(f) = \int_{\mathbb{Z}_p} f(x)   d\mu_{-q} (x)= \lim_{N \rightarrow \infty} \frac{1}{[p^N]_{-q}} \sum_{x=0}^{p^N-1} f(x) (-q)^x.
\end{split}\end{equation}
From \eqref{32}, we note that
\begin{equation}\begin{split}\label{33}
qI_{-q}(f_1) +I_{-q}(f) = [2]_q f(0),\,\,\text{where}\,\, f_1(x) = f(x+1).
\end{split}\end{equation}
By \eqref{33}, we get
\begin{equation}\begin{split}\label{34}
\left( \frac{q+(1+\lambda t)^{-\frac{1+q}{\lambda } }}{q} \right)\int_{\mathbb{Z}_p}   (1+\lambda t)^{-\frac{x}{\lambda }(1+q)} d\mu_{-q^{-1}} (x) = [2]_{q^{-1}},
\end{split}\end{equation}
where $\lambda \in \mathbb{Z}_p$ and $|t|_p < p^{-\frac{1}{p-1}}$. Thus, from \eqref{34}, we have
\begin{equation}\begin{split}\label{35}
\int_{\mathbb{Z}_p} (1+\lambda t)^{-\frac{x}{\lambda }(1+q)}   d\mu_{-q^{-1}} (x) = \frac{1+q}{q+(1+\lambda t)^{-\frac{1+q}{\lambda }}}
\end{split}\end{equation}
From \eqref{16} and \eqref{35}, we note that
\begin{equation}\begin{split}\label{36}
\int_{\mathbb{Z}_p} (1+\lambda t)^{-\frac{x}{\lambda }(1+q)}   d\mu_{-q^{-1}} (x)= \sum_{n=0}^\infty A_{n,\lambda }(-q)\frac{t^n}{n!}.
\end{split}\end{equation}
Now, we define the degenerate rising factorials as follows:
\begin{equation}\begin{split}\label{37}
<x>_{0,\lambda }=1, \,\,<x>_{n,\lambda } = x(x+\lambda )(x+2\lambda )\cdots (x+(n-1)\lambda ),\,\,(n \geq 1).
\end{split}\end{equation}
It is not difficult to show that
\begin{equation}\begin{split}\label{38}
(1+\lambda t)^{-\frac{x}{\lambda }(1+q)}&= e^{-\frac{x}{\lambda }(1+q)\log(1+\lambda t)} = \sum_{n=0}^\infty (-1)^n \left( \frac{x}{\lambda } \right)^n (1+q)^n \frac{\big(\log(1+\lambda t)\big)^n}{n!}\\
&=\sum_{k=0}^\infty (-1)^k \left( \frac{x}{\lambda } \right)^k (1+q)^k \sum_{n=k}^\infty S_1(n,k) \lambda ^n \frac{t^n}{n!}\\
&= \sum_{n=0}^\infty \left( \sum_{k=0}^n (-1)^k \lambda ^{n-k} (1+q)^k S_1(n,k) x^k \right)  \frac{t^n}{n!}.
\end{split}\end{equation}
From \eqref{33}, we note that
\begin{equation}\begin{split}\label{39}
\int_{\mathbb{Z}_p}  e^{xt}  d\mu_{-q^{-1}} (x) = \frac{q+1}{e^t+q} = \sum_{n=0}^\infty H_n(-q) \frac{t^n}{n!}.
\end{split}\end{equation}
Thus, by \eqref{39}, we get
\begin{equation*}\begin{split}
\int_{\mathbb{Z}_p} x^n    d\mu_{-q^{-1}} (x) = H_n(-q),\,\,(n \geq 0),
\end{split}\end{equation*}
where $H_n(-q)$ are the Frobenius-Euler numbers. From \eqref{38} and \eqref{39}, we have
\begin{equation}\begin{split}\label{40}
\int_{\mathbb{Z}_p} (1+\lambda t)^{-\frac{x}{\lambda }(1+q)}   d\mu_{-q^{-1}} (x) = \sum_{n=0}^\infty \left( \sum_{k=0}^n (-1)^k \lambda ^{n-k} (1+q)^k S_1(n,k) H_k(-q) \right) \frac{t^k}{k!}
\end{split}\end{equation}
Comparing the coefficients on both sides of \eqref{36} and \eqref{40}, we have
\begin{equation}\begin{split}\label{41}
A_{n,\lambda }(-q)= \sum_{k=0}^n (-1)^k \lambda ^{n-k} (1+q)^k S_1(n,k) H_k(-q),\,\,(n \geq 0).
\end{split}\end{equation}
In particular,
\begin{equation}\begin{split}\label{42}
(1+\lambda t)^{-\frac{x}{\lambda }(1+q)} &= \sum_{n=0}^\infty {-\frac{x}{\lambda }(1+q) \choose n} \lambda ^n t^n\\
&= \sum_{n=0}^\infty \big(-\tfrac{x}{\lambda }(1+q)\big)_n\lambda ^n \frac{t^n}{n!}\\
&= \sum_{n=0}^\infty (-1)^n <(1+q)x>_{n,\lambda } \frac{t^n}{n!}\\
&= \sum_{n=0}^\infty (-1)^n <x>_{n,\frac{\lambda }{1+q}} (1+q)^n \frac{t^n}{n!}
\end{split}\end{equation}
From \eqref{42}, we note that
\begin{equation}\begin{split}\label{43}
\sum_{n=0}^\infty A_{n,\lambda }(-q) \frac{t^n}{n!} &= \int_{\mathbb{Z}_p}  (1+\lambda t)^{-\frac{x}{\lambda }(1+q)}   d\mu_{-q^{-1}} (x)\\
&= \sum_{n=0}^\infty (-1)^n (1+q)^n \int_{\mathbb{Z}_p} <x>_{n,\frac{\lambda }{1+q}}   d\mu_{-q^{-1}} (x)\frac{t^n}{n!}.
\end{split}\end{equation}
Thus, by comparing the coefficients on the both sides of \eqref{43}, we get
\begin{equation}\begin{split}\label{44}
\int_{\mathbb{Z}_p}  <x>_{n,\frac{\lambda }{1+q}}  d\mu_{-q^{-1}} (x)= (-1)^n \frac{A_{n,\lambda }(-q)}{(1+q)^n}, \,\,(n \geq 0).
\end{split}\end{equation}
For any positive real number  $\lambda$, the degenerate unsigned Stirling numbers of the first kind $| S_{1,\lambda }(n,l) |$ are defined by
\begin{equation}\begin{split}\label{45}
<x>_{n,\lambda }= \sum_{l=0}^n | S_{1,\lambda }(n,l) | x^l,\,\,(n \geq 0).
\end{split}\end{equation}
From \eqref{45}, we have
\begin{equation}\begin{split}\label{46}
\int_{\mathbb{Z}_p} <x>_{n,\frac{\lambda }{1+q}}  d\mu_{-q^{-1}} (x)&= \sum_{l=0}^n |S_{1,\frac{\lambda }{1+q}} (n,l)| \int_{\mathbb{Z}_p} x^l   d\mu_{-q^{-1}} (x)\\
&= \sum_{l=0}^n |S_{1,\frac{\lambda }{1+q}} (n,l)|H_l(-q).
\end{split}\end{equation}
Hence, by \eqref{44} and \eqref{46}, we get
\begin{equation*}\begin{split}
\sum_{l=0}^n |S_{1,\frac{\lambda }{1+q}} (n,l)|H_l(-q) = (-1)^n \frac{A_{n,\lambda }(-q)}{(1+q)^n},\,\,(n \geq 0).
\end{split}\end{equation*}

\end{document}